\documentclass[final]{aupl}

\usepackage{
amsfonts,
latexsym,
amssymb,
enumerate,
verbatim,
mathrsfs,
graphicx,
centernot,
}

\newcommand{\labbel}{\label}

\newtheorem{theorem}{Theorem}[section]
\newtheorem{lemma}[theorem]{Lemma}
\newtheorem{thm}[theorem]{Theorem}

\newtheorem{corollary}[theorem]{Corollary}

\newtheorem*{claim*}{Claim}

\newtheorem*{theorem*}{Theorem}
\newtheorem*{proposition*}{Proposition}
\newtheorem*{corollary*}{Corollary}
\newtheorem*{lemma*}{Lemma}
\newtheorem*{scholion*}{Scholion}

\theoremstyle{definition}
\newtheorem{definition}[theorem]{Definition}

\newtheorem{problem}[theorem]{Problem}

\theoremstyle{remark}
\newtheorem{remark}[theorem]{Remark}
\newtheorem{remarks}[theorem]{Remarks}
\newtheorem*{remark*}{Remark}
\newtheorem*{remarks*}{Remarks}

\newtheorem*{observation*}{Observation}

 \allowdisplaybreaks[1]

\numberwithin{equation}{section}

\begin{document}

\title[Linear orders with one operator,
amalgamation property]{Linearly ordered sets with only one operator
have the amalgamation property}

\author{Paolo Lipparini} 
\address{Dipartimento di Matematica\\Viale della  Ricerca
 Scientifica Lineare\\Universit\`a di Roma ``Tor Vergata'' 
\\I-00133 ROME ITALY}

\email{lipparin@axp.mat.uniroma2.it}



\date{\today}

\begin{abstract}
The class of linearly ordered sets with one
order preserving unary operation 
has the Strong Amalgamation Property (SAP).
The class of linearly ordered sets with one strict order preserving
  unary operation
has AP but not SAP.
The class of linearly ordered sets with 
 two  order preserving  unary operations
has not AP.
For every set $F$, the class of linearly ordered sets
with an $F$-indexed family of automorphisms has SAP.   
Corresponding results are proved in the case of order reversing operations.
Various subclasses of the above classes  are considered
and some model-theoretical consequences are presented.
 \end{abstract}

\maketitle

\section{Introduction} \labbel{intro}

The amalgamation property (AP) has found deep 
applications in algebra and logic,
and is nontrivially linked to categorical notions.
In the special case of groups 
the amalgamation property has been 
considered in Schreier \cite{S}. 
 Then Fra\"\i ss\'e \cite{F1,F}
and J{\'o}nsson \cite{J}  
introduced the abstract general definitions 
and initiated a flourishing line of research
with applications in model theory.
Subsequently, another  line of research 
connected the amalgamation property with algebraic logic.
See \cite{E,GM,H,KMPT,MMT} for 
more details and further references.

If one
starts with some theory having AP
and adds a set of 
operators with suitable properties,
 sometimes the resulting theory has still AP.
Many  results of this kind are known for fields
with operators, e.~g.,
\cite{W,Z}.
See \cite{BHKK,GP} for more recent results and further references. 
A similar preservation phenomenon sometimes occurs for
ordered structures \cite{LP}.
In particular, many kinds of Boolean algebras 
with operators have the strong amalgamation 
property (SAP). See 
\cite{J2,EC,MS,N}\footnote{Results in \cite[Theorem 4]{EC}
are stated for just one operator, but the proof
works for an arbitrary set of operators.}.
The case of partially ordered sets
with any number of order preserving unary operations is probably folklore;
anyway, see Corollary \ref{oper} below.

In Section \ref{lop} we prove
the quite curious fact that, 
on the other hand, for \emph{linearly} ordered sets,
adding a single order preserving unary operation maintains SAP, 
while AP fails  when two operations are present.
Linearly ordered sets with one strict order
preserving operation have AP but not SAP.
Corresponding results are proved in Section \ref{r} for 
linearly ordered sets with order reversing  operations,
but in this case  strong  
amalgamation  generally fails.
On the positive side, any class
of linearly ordered sets with  families of automorphisms
and antiautomorphisms with a common fixed point
has SAP.

An appealing aspect of our proofs 
is that  we always construct the amalgamating structure on 
the set-theoretical union of the domains of the structures
to be amalgamated. Henceforth we need no
further effort in order to  get AP  for various subclasses
of the above classes. The existence of Fra\"\i ss\'e 
limits follows in most cases, and sometimes we even
get model completions for appropriate theories.
This aspect is discussed in Section \ref{fur}.

 \subsection{Outline of the proofs} \labbel{out}  
Our main techniques are summarized as follows.
Given linearly ordered sets 
$\mathbf A, \mathbf B $ and $  \mathbf C $
to be amalgamated, with ground $\mathbf C$,
we first use  Fra\"\i ss\'e \cite{F}
and J{\'o}nsson's \cite{J} method in order to embed
$\mathbf A $ and $  \mathbf B$
into a partially ordered set $\mathbf E$ over $A \cup B$. 
In the absence of operations, it is enough to extend
the partial order on $\mathbf E$ to some arbitrary linear order,
but a cleaner (classical) way to do this is to consider
some element of $A$ to be always smaller than  
some element of $B$, provided there is no other relation
to be satisfied and which implies the converse.

Given a triple to be amalgamated, 
the above idea provides a rather uniform method to extend 
a partial order, and the method works even in the presence of
one operation $f$.
Of course, we are not allowed to always set
 $a < b$, whenever $a \in A$ and $b \in B$
are not comparable in $\mathbf E$.
Indeed it may happen that  $f(a)$ and $f(b)$
  are  comparable in $\mathbf E$
and in this case the relative position of $a$
and $b$ should be set accordingly.  
We check that all the  conditions arising in a similar way
can be consistently put together, hence we succeed
in getting a linear order.

The most delicate
case is when the operation is supposed
to be strict order preserving.
In this case, some elements of 
$A$ and $B$ possibly need to be identified:
this means that \emph{strong} amalgamation
 fails. However, the relations
involved in the identifications exactly determine
the structure relative to such particular elements.
In other words, the ground structure
can be extended to some model $\mathbf C_1$
which then becomes a strong amalgamation base. 
Needless to say, details are delicate  in each case,
since the method works only for one operation,
not for a pair of operations, the counterexamples
being quite easy.
In contrast, and  to make the situation even more involved,
SAP holds for any number of operations,
under the assumption that the operations are 
automorphisms.

The arguments need to be modified when dealing with order
reversing operations. In this case, we cannot
always put $a < b$, when $a \in A$, $b \in B$ 
and the ordering relation between $a$ and $b$ is not determined by
other conditions. Indeed, if $g$
is order reversing, then $a < b$ implies $g(b) \leq g(a)$, 
but still  $g(a) \in A$ and  $g(b) \in B$.
However, in the presence of an order reversing operation,
the elements of a linearly ordered set can be obviously divided into
 ``lower''   and ``upper''  elements: we set 
$a < b$ for lower elements and $b < a$
for upper elements, again, when there is no other condition
to be satisfied.

\section{Preliminaries} \labbel{prel}

In this section we recall the basic definitions and
some classical constructions which show amalgamation
for partial and linear orders.

\begin{definition} \labbel{sap}
If $\mathcal K$ is a
class of structures  of the same type,
then $\mathcal K$ is said to have the \emph{amalgamation property}
(AP) if, whenever 
$\mathbf A, \mathbf B, \mathbf C \in \mathcal K$,
 $ \iota _{\mathbf C, \mathbf A} \colon \mathbf C \rightarrowtail \mathbf A$
and 
 $ \iota _{\mathbf C, \mathbf B} \colon \mathbf C \rightarrowtail \mathbf B$
are embeddings, then there is a structure
$\mathbf D \in \mathcal K$ and  embeddings
$ \iota _{\mathbf A, \mathbf D} \colon \mathbf A \rightarrowtail \mathbf D$
and 
 $ \iota _{\mathbf B, \mathbf D} \colon \mathbf B \rightarrowtail \mathbf D$
such that 
$\iota _{\mathbf C, \mathbf A} \circ \iota _{\mathbf A, \mathbf D}
=
\iota _{\mathbf C, \mathbf B} \circ \iota _{\mathbf B, \mathbf D}$. 
Namely, the following diagram can be commutatively completed as requested.
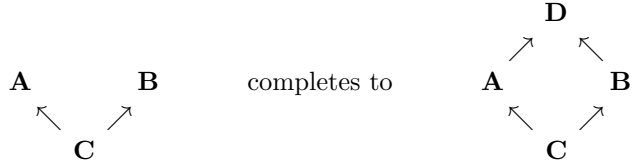
\begin{equation*}
\begin{matrix} 
 \cr
 \cr
\mathbf A \quad \quad \quad \quad \mathbf B \cr
\nwarrow \ \quad\  \nearrow \cr
 \mathbf C
 \end{matrix}   
\qquad \quad  \text{ completes to }  \quad \qquad 
\begin{matrix} 
 \mathbf D \cr
\nearrow \ \quad\  \nwarrow \cr
\mathbf A \quad \quad \quad \quad \mathbf B \cr
\nwarrow \ \quad\  \nearrow \cr
 \mathbf C
 \end{matrix}   
\end{equation*}

 If, in addition, 
the above model and embeddings can be
always chosen in such a way that
 the intersection
of the images of $\iota _{\mathbf A, \mathbf D}$
and $\iota _{\mathbf B, \mathbf D}$ 
is equal to the image of 
$\iota _{\mathbf C, \mathbf A} \circ \iota _{\mathbf A, \mathbf D}$, then
$\mathcal K$ is said to have the \emph{strong amalgamation property}
(SAP).

The latter condition
can be simplified under the assumption that 
$\mathcal K$ is  closed under isomorphism.
Under this assumption, $\mathcal K$ 
has SAP  if and only if, whenever 
$\mathbf A, \mathbf B, \mathbf C \in \mathcal K$,
 $ \mathbf C \subseteq  \mathbf A$,
 $ \mathbf C \subseteq  \mathbf B$
and $C=A \cap B$,
then there is a structure
$\mathbf D \in \mathcal K$ such that 
$ \mathbf A \subseteq  \mathbf D$ and 
$ \mathbf B \subseteq  \mathbf D$.
Here, say, $ \mathbf C \subseteq  \mathbf A$ 
means that $C \subseteq A$ as sets and that the inclusion
is an embedding   from $ \mathbf C $ to $  \mathbf A$. 

A triple of models  $ \mathbf A $, $  \mathbf B$
and  $ \mathbf C $ as above shall be called an
\emph{amalgamation triple}, or a \emph{triple to be amalgamated}.
The structure  $ \mathbf  D$ shall be called a \emph{(strong) amalgam}, or
a \emph{(strong) amalgamating structure}.  
Whenever possible, we shall consider the simplified setting
described in the previous paragraph, namely, we shall deal with
inclusions rather than with arbitrary embeddings.
The setting in which we  work
 shall always be clear from the context.

Even if some class $\mathcal K$ has not AP,
it is anyway interesting to ask when a diagram as above
can be completed. 
In particular, it is interesting to consider those
specific  $ \mathbf C $
for which the diagram can be always completed.
In detail, a structure 
$ \mathbf C $ 
is said to be a \emph{(strong) amalgamation base 
for a class $\mathcal K$} 
if  every amalgamation triple with 
 $ \mathbf C $ at the bottom and $\mathbf A, \mathbf  B \in \mathcal K$ 
 has some (strong) amalgamating structure $\mathbf  D \in \mathcal K$.

Notice that here we are always dealing with \emph{embeddings},
not \emph{homomorphisms}.
E.~g., in the case of ordered sets,
 $ \iota \colon \mathbf C \rightarrowtail \mathbf A$
is an embedding if, for every $c,d \in C$,
it happens that $c \leq _{\mathbf C} d$
if and only if 
 $\iota(c) \leq _{\mathbf A} \iota(d)$.
In the definition of an homomorphism the sole
``only if'' implication is required.
Homomorphisms of partially ordered sets are frequently called
\emph{ordermorphisms}.

In the above inequalities we have written
$ \leq _{\mathbf C} $ 
and $ \leq _{\mathbf A} $ 
to distinguish the order relation considered on $\mathbf C$ 
from the order relation considered on $ \mathbf A$.
We shall use a similar convention when dealing with
 unary operations. 
As customary, we shall drop the subscripts when there is
no risk of confusion.

As a final detail, slightly distinct
notions arise if in (S)AP
one allows or does not allow $\mathbf C$ to be an empty structure. 
The results here shall not be affected by the distinction,
hence the reader might use her or his favorite
version of the definition.
\end{definition}

The following classical construction
will be the starting point for our proofs.
\emph{Poset}
is an abbreviation for \emph{partially ordered set}.

\begin{theorem} \labbel{FJ}
(Fra\"\i ss\'e \cite[9.3]{F}, J{\'o}nsson  \cite[Lemma 3.3]{J}) 
The class 
%$\mathcal{PO} $
 of posets
has the strong amalgamation property.

If $\mathbf A, \mathbf B, \mathbf C $ are posets,
 $ \mathbf C \subseteq  \mathbf A$,
 $ \mathbf C \subseteq  \mathbf B$
and $C=A \cap B$,
then an amalgamating structure
$\mathbf E $  is obtained as follows.
The domain of $\mathbf E $ is $A \cup B$ and,
for $d,e \in A \cup B$, let $d \leq _{\mathbf E} e $
if either   
\begin{gather}
\begin{aligned} \labbel{al}   
& d,e \in A \text{ and } d \leq_{\mathbf A} e, \text{ or }
\\  
& d,e\in B \text{ and } d \leq_{\mathbf B} e, \text{ or }
\\  
& d \in A, e \in B \text{ and 
there is } c \in C  \emph{ such that }
 d \leq_{\mathbf A} c \text{ and } c \leq_{\mathbf B} e, \text{ or }
\\  
& d \in B, e \in A \text{ and 
there is } c \in C  \emph{ such that }
 d \leq_{\mathbf B} c \text{ and } c \leq_{\mathbf A} e.
\end{aligned}
  \end{gather}   

For short, 
${\leq _{\mathbf E}}=
{\leq _{\mathbf A}} \cup {\leq _{\mathbf B}} \cup
({\leq _{\mathbf A}} \circ {\leq _{\mathbf B}})
\cup  ({\leq _{\mathbf B}} \circ {\leq _{\mathbf A}})$.
 \end{theorem} 

See \cite{F,J} for a proof and
\cite{apu} for various generalizations.  

The order $\leq _{\mathbf E}$
as defined above is the finest,
i.e., smallest  order
on $D= A \cup B$ which makes $D$
an amalgamating structure. However, 
$\leq _{\mathbf E}$ is not the unique such order.
For example, if $\mathbf A, \mathbf B $ and $  \mathbf C $
 are linear orders, then, for 
any extension $\leq _{\mathbf D}$ of  $\leq _{\mathbf E}$
to a linear order, the inclusions
from $\mathbf A$ and $\mathbf  B$ to 
$(A \cup B, \leq _{\mathbf D})$ 
 still remain embeddings.
This is enough to show that the class of linearly ordered sets
has SAP, but a cleaner method is to extend  $\leq _{\mathbf E}$ 
in such a way that
$a<_{\mathbf D} b$, whenever  $a \in A$, $b \in B$ 
and the relative order between $a$ and  $b$
is not decided by $\leq _{\mathbf E}$.
This is a classical argument, see, e.~g., \cite[9.2]{F},
\cite[Example 2.2.1]{E}. We present full details in the   
next corollary, since similar methods will prove useful
in the following sections.

\begin{corollary} \labbel{losp}
The class 
%$\mathcal {LO}$
 of linearly ordered sets has SAP.
 \end{corollary} 

 \begin{proof}
Let $\mathbf C \subseteq \mathbf A, \mathbf B$ 
be a triple of linear orders to be amalgamated.
In particular, 
$\mathbf A$, $\mathbf B$ and $\mathbf C$
are partial orders, hence Theorem \ref{FJ} 
provides an amalgamating \emph{partial} order
$\mathbf E=(D, \leq_{\mathbf E} )$ with
$D=A \cup B$.
If $d, e \in D$ are $\leq_{\mathbf E}$-incomparable, 
then    necessarily 
$d \in A$ and $e \in B$, or conversely, since
$\mathbf A$ and  $\mathbf B$ are linearly ordered
and $D= A \cup B$.
Extend $\leq_{\mathbf E}$ to a relation 
$\leq_{\mathbf D}$ by always letting the element in $A$
be  $<_{\mathbf D}$ than the element in $B$, for every 
pair of $\leq_{\mathbf E}$-incomparable elements. 
In view of \eqref{al}, it is easy to check
transitivity and antisymmetry of $\leq_{\mathbf D}$.
The resulting order is linear by construction.

Letting    $\mathbf D=(D, \leq_{\mathbf D})$,
the identity map $\iota$  is an ordermorphism
from $\mathbf E$ to $\mathbf D$, though  not necessarily an
embedding. However, the composition   
$\iota _{\mathbf A, \mathbf E} \circ \iota $
is indeed an embedding from $\mathbf A$ to $\mathbf D$,
 and similarly for $\mathbf B$,
hence   $\mathbf D$ is an amalgamating
structure in the class of \emph{linear} orders. 

As we mentioned, we could have done by simply
extending the partial order $\leq_{\mathbf E}$ 
to some arbitrary linear order $\leq_{\mathbf F}$;
however, the assumption that every partial order
can be extended to a linear order is a weak
form of the axiom of choice \cite[Form 49]{HR};
on the other hand, the  argument we have recalled
seems to need no form of choice.
Moreover, the above method, with variations, shall be used 
in order to prove Theorems \ref{los} and \ref{rev}  below.  
 \end{proof} 

As another immediate consequence of Theorem \ref{FJ},
if we add monotone unary operations to  partial orders,
then SAP is maintained. See the next corollary. We shall see in the following sections
that this is not always the case, when dealing with linear orders.
 
If $\mathbf A$ is a poset,  a unary operation
$f \colon A \to A$ is \emph{order preserving}, resp., 
\emph{order reversing}
if, for every $a, b \in A$,  $a \leq b $ implies
$f(a) \leq f(b)$, resp., $f(a) \geq f(b)$.
We say that $f$ is
 \emph{strict order preserving}, resp., 
\emph{strict order reversing}
if, for every $a, b \in A$,  $a < b $ implies
$f(a) < f(b)$, resp., $f(a) > f(b)$.

\begin{corollary} \labbel{oper}
The class  of posets with any 
(fixed in advance) number of 
order preserving, order reversing, 
strict order preserving and strict order reversing
unary operations has SAP, actually, the superamalgamation 
property; see \cite[p. 173]{GM}.
 \end{corollary} 

More formally, Corollary \ref{oper}
asserts that, for every set 
$F= F_1 \cup F_2 \cup F_3 \cup F_4$  
of unary function symbols, 
if $\mathcal {PO}_F$ is 
the class of posets 
with additional functions such that 
all the symbols in $F_1$ are interpreted
as order preserving functions, all the symbols in $F_2$ are interpreted
as order reversing functions, etc.,
 then  $\mathcal {PO}_F$
 has SAP.

\begin{proof}
Given three structures 
$\mathbf A, \mathbf B, \mathbf C $  
for the appropriate language
and to be amalgamated,
first construct a partial order $\leq _{\mathbf E} $ 
on $D=A \cup B$ 
as in Theorem \ref{FJ}.
Since  $D=A \cup B$ and the operations
under consideration are unary and agree on $C=A \cap B$, then
each operation can be uniquely extended over $D$.
Notice that here it is fundamental to have the strong
version of the amalgamation property for posets.
It is immediate from 
\eqref{al} that if, say, 
$f$ is interpreted by an order preserving
operation on $\mathbf A, \mathbf B, \mathbf C $,
 then the extension of $f _{\mathbf A} $
and $f _{\mathbf B} $ to $\mathbf E$
is still order preserving. 
For example, if $a \leq_{\mathbf E} b$
is given by $a \leq_{\mathbf A} c \leq_{\mathbf B} b$, for some
$c \in C$, then    
$f_{\mathbf A}(a) \leq_{\mathbf A} f_{\mathbf A}(c)$ and 
$f_{\mathbf B}(c) \leq_{\mathbf B} f_{\mathbf B}(b)$,
since $f_{\mathbf A}$ and  $f_{\mathbf B}$
are order preserving in  ${\mathbf A}$ and  ${\mathbf B}$,
respectively. Hence $f(a) \leq_{\mathbf E} f(b)$.
 \end{proof}

\section{Linearly ordered sets with operators} \labbel{lop} 

We have seen in the previous
section that posets, possibly with 
operators, always share SAP.
 The situation with linearly ordered sets is 
much more delicate.
Everything runs smoothly
when at most one order preserving operation is added,
but even AP fails when more than one operation
are present.
One strict order preserving operation
prevents SAP, but AP is still satisfied.
In contrast, SAP holds when an arbitrary family of automorphisms
is added.

The method recalled in the proof of 
Corollary \ref{losp} needs to  be modified, since
it is possible that, say, 
$a$ and  $b$ are not comparable in $\mathbf E$,
while
 $f(a)$ and $f(b)$
  turn out to be comparable.
We check that all the  conditions arising in similar ways
can be consistently put together in the case of just an operation, 
while this is not possible for two or more operations.

As usual, if $f$ is a unary operation, 
$f^n$ is defined inductively 
as follows: $f^0(a)= a$
and  $f^{n+1}(a)= f(f^n(a))$.

\begin{theorem} \labbel{los} 
  \begin{enumerate}[(a)]
\item
 The class 
$\mathcal {LO}_p$
 of linearly ordered sets with an order preserving
unary operation 
has SAP.

\item The class 
$\mathcal {LO}_{sp}$ 
of linearly ordered sets with a strict order preserving
unary operation 
has AP but not SAP.

\item
The classes 
%$\mathcal {LO}_{2p}$, resp., $\mathcal {LO}_{2sp}$ 
of linearly ordered sets with 
 two order preserving, 
resp.,  two strict order preserving
unary operations
have not AP. 

\item
For every set $ F$, the class 
$\mathcal {LO}_{Fa}$
 of 
linearly ordered sets with an $F$-indexed
family of  order automorphisms
has SAP. 
  \end{enumerate} 
\end{theorem} 

 \begin{proof}
(a)  Fix
a triple $\mathbf C \subseteq \mathbf A, \mathbf B$
to be amalgamated.
The proofs of 
 Theorem \ref{FJ} and Corollary  \ref{oper}  
furnish an amalgamating  \emph{partial} order
$\mathbf E=(D, \leq _ {\mathbf E} , f)$ on 
$D=A \cup B$ such that 
$f$ is an  order preserving operation  with respect to 
$\leq _ {\mathbf E}$.
Szigeti and Nagy \cite{SN} 
provided a condition under which a partial order $D$ 
with an order preserving unary operation $f$ can
be linearized in such a way that the operation
is still order preserving with respect to the new linear 
order. This holds if and only if 
$f$ is \emph{acyclic}, namely,
whenever $d \in D$ and  
$f ^{m+1}(d) \in \{d, f (d), f ^{2}(d), \dots, f ^{m}(d)  \}  $,
for some $m \in \mathbb N$, then    
$f ^{m+1}(d)=f ^{m}(d)$. 
In the case at hand, $f$ is surely acyclic, 
since $D=A \cup B$ and $f$   is 
trivially acyclic on both $A$ and $B$.
Relying on \cite{SN} is thus sufficient in order
to prove (a). To make the paper self-contained, we shall also present a more direct 
proof of (a) which has the advantage of making no use
of the axiom of choice.

So let again 
$\mathbf E=(D, \leq _ {\mathbf E} , f)$ be given by
 Theorem \ref{FJ} and Corollary  \ref{oper}, with  
$D=A \cup B$ and
$\leq _ {\mathbf E}$
only a partial order.
As in the proof
of Corollary  \ref{losp}, we shall extend 
 $\leq _ {\mathbf E}$ to some linear order $\leq$ on $D$,
but in the present case the behavior of the operation 
$f$ should be taken into account. 
The values of $f$  shall not be modified.

We first explicitly describe the way  the linear order on 
$\mathbf C$ forms the ``backbone''
of the partial order $\leq _ {\mathbf E}$.
This description does not involve $f$
and in principle is not strictly necessary,
but it will greatly simplify the subsequent arguments. 
Recall that if $\mathbf C$
is a linearly ordered set, a \emph{cut} of  $\mathbf C$
is a pair $(C_1, C_2)$ such that 
$C_1 \cup C_2 = C$
and $c_1 < c_2$, for every     
$c_1 \in C_1$ and $c_2 \in  C_2$,
in particular, $C_1 \cap C_2 = \emptyset $.
We allow $C_1$ or $C_2$ to be empty.  
If $\mathcal G = (C_1, C_2)$ is a cut  of $\mathbf C$,
 the \emph{component (of $D$)
associated to}
$\mathcal G$ is the set 
$O_{\mathcal G} = \{ \, d \in D  \mid c_1 < d < c_2, \text{ for all } c_1 \in C_1
\text{ and }   c_2 \in C_2 \,\}$.
Thus  $O_{\mathcal G} \subseteq D \setminus C$,
since  $C=C_1 \cup C_2 $.
Actually, every element  $d$ of $D \setminus C$
belongs  to some component:
if $d \in D \setminus C$, then  
$d$ \emph{determines} a cut $\mathcal G $ of $\mathbf C$
by setting 
$C_1 = \{ \,  c \in C \mid c < _ {\mathbf E} d \, \} $
and 
$C_2 = \{ \,  c \in C \mid d < _ {\mathbf E} c \, \} $.
Then $\mathcal G = (C_1, C_2)$ is a cut, since $d \notin C$,
$\mathbf C$ embeds in $\mathbf E$ 
%%%CAMBIATO QUI
and 
$ \mathbf A$,  $ \mathbf B$ are linearly ordered.  
If $d$ determines $\mathcal G $, then $d$ belongs to the component 
associated to $\mathcal G $; actually, this is the only component
to which $d$ belongs. 
In other words, the nonempty components partition 
$ D \setminus C$.

If $d$ and $d'$ lie in two distinct components, then
either $ d \leq _ {\mathbf E} d'$
or $d' \leq _ {\mathbf E} d$.
Indeed, let   
$(C_1, C_2)$ and
$(C'_1, C'_2)$
be the cuts determined by 
$d$ and $d'$, thus
$C_1 \neq C'_1$.
Since
$C_1 \cup C_2 = C'_1 \cup C'_2 = C$, then either 
$C_1 \cap C'_2 \not= \emptyset $ or 
$C'_1 \cap C_2 \not= \emptyset $.
If $c \in C_1 \cap C'_2$, then
$d' \leq _ {\mathbf E} c \leq _ {\mathbf E} d$.
Similarly, if $c \in C'_1 \cap C_2$, then
$d \leq _ {\mathbf E} c \leq _ {\mathbf E} d'$.

Moreover, if $c \in C$ and $d \in D$, then either  
$ c \leq _ {\mathbf E} d$
or $d \leq _ {\mathbf E} c$,
since 
$\leq _ {\mathbf E}$ extends the linear orders
$\leq _ {\mathbf A}$ and $\leq _ {\mathbf B}$,
$D= A \cup B$, $c \in C = A \cap B$, hence
either $c,d \in A$ or $c,d \in B$, thus
 the relative position of $c$
and $d$ is already determined by
either $\leq _ {\mathbf A}$ or $\leq _ {\mathbf B}$.  
It follows that 

(*) in order to extend $\leq _ {\mathbf E}$
to a linear order on $D$ it is enough
to extend to a linear order each  restriction of 
$\leq _ {\mathbf E}$  to each component.

So let $O$ be the component associated to some cut.
Recall that $O \subseteq D \setminus C$.  
If $a,a' \in O \cap A$, then either  
$ a \leq _ {\mathbf E} a'$
or $ a' \leq _ {\mathbf E} a$, by \eqref{al}, since
$\leq _ {\mathbf A}$ is a linear order,
$\leq _ {\mathbf E}$  coincides with $\leq _ {\mathbf A}$ 
on $A$  and $a$, $a'$ lie in the same component.     
The situation is similar if $b,b' \in O \cap B$.
If $a \in O \cap A$
and  $b \in O \cap B$,
then $a$ and  $b$
are  $  \leq _ {\mathbf E} $-incomparable, by the last two lines in
condition \eqref{al} and since $a$ and $b$ lie in the same component.   
Henceforth it is enough to set the relative order 
for each pair $a \in A \cap O$ and $b \in B \cap O$.
  
Let $<_O$ extend  $  < _ {\mathbf E} $ on $O$ by setting
\begin{align} \labbel{alal}   
 &\left. a<_O a_1,\right. 
&& \text{if $a, a_1 \in O \cap A$ and
$a <_ {\mathbf A} a_1$}
\\
\labbel{alalu}
&\left.b<_O b_1,\right. 
&& \text{if $b,b_1 \in O \cap B$ and  
$b < _ {\mathbf B} b_1$}
\\
\labbel{alaluu} 
&\left. \begin{aligned}  
 &b <_O a, 
&&\text{if $ f^n(b) <_ {\mathbf E} f^n(a)  $},
\\[- 2pt]&&&\text{for some
$n \geq 1$, or} 
\\
& a<_O b, &&\text{otherwise,}
\end{aligned} 
\ \ \right\}
&&\text{if $a \in A \cap O$ and $b \in B \cap O$}
  \end{align}

We now show that 
$<_O$ is a linear order on $O$.  
As above, all pairs of distinct elements 
in $A \cap O$ are $  < _ {\mathbf A} $ comparable, 
hence $  < _ {O} $ comparable. 
A similar remark holds for 
$B \cap O$. By construction, all pairs
$a \in A \cap O$ and $b \in B \cap O$
are $<_O$ comparable and we cannot have both
$a <_O b$ and  $b <_O a$.
Notice that $A \cap B \cap O = \emptyset $,
since $A \cap B = C$
and  $O$ is a component, hence $O \subseteq D \setminus C$.    
Moreover, $d <_O d$ is impossible, if $d \in O$. 

It remains to show that 
$<_O$ is transitive on $O$.
The proof goes by considering all possible cases.
Notice that Clause \eqref{alaluu}  is not symmetrical. Let 
$a, a_1 \in A \cap O$ and $b, b_1 \in B \cap O$.

(i)
If $a <_O a_1 <_O b$, then
for no $n$ $ f^n(b) <_ {\mathbf E} f^n(a_1)  $,
a fortiori   for no $n$ $ f^n(b) <_ {\mathbf E} f^n(a)  $.
Indeed, $a <_O a_1$ means
$ a <_ {\mathbf A} a_1  $
hence 
$ f^n(a) \leq _ {\mathbf A} f^n(a_1)  $,
since $f$ is order preserving on $\mathbf A$, thus 
$ f^n(a) \leq _ {\mathbf E} f^n(a_1)  $.
If by contradiction
$ f^n(b) <_ {\mathbf E} f^n(a)  $, then
$ f^n(b) <_ {\mathbf E} f^n(a) \leq _ {\mathbf E} f^n(a_1) $,
 a contradiction, since 
$ \leq _ {\mathbf E}$ is a partial order.
Hence $a <_O b$.

(ii)
If $a <_O b <_O a_1$, then
$ f^n(b) <_ {\mathbf E} f^n(a_1)  $, for some
$n $.
Were $a_1 \leq_O a$, then 
$ f^n(a_1)  \leq _ {\mathbf A} f^n(a)  $,
since $f$ is order preserving on $\mathbf A$, 
thus 
$ f^n(b) <_ {\mathbf E}  f^n(a_1)  \leq _ {\mathbf E} f^n(a)  $,
contradicting $a <_O b $.
Hence $a <_O   a_1$.
Here we have used the assumption that
$ \leq_ {\mathbf A} $ is a linear order on $A$.

(iii)
If $a <_O b <_O b_1$, then
for no $n$ $ f^n(b) <_ {\mathbf E} f^n(a)  $,
a fortiori   for no $n$ $ f^n(b_1) <_ {\mathbf E} f^n(a)  $,
since $ f^n(b) \leq_ {\mathbf E} f^n(b_1)  $.
Hence $a  <_O b_1$.

(iv)
If $b <_O a <_O a_1$, then
$ f^n(b) <_ {\mathbf E} f^n(a)  $, for some
$n $, hence $ f^n(b) <_ {\mathbf E} f^n(a) \leq_ {\mathbf E} f^n(a_1)  $,
thus 
$b <_O a_1$.

(v)
If $b <_O a <_O b_1$, then
$ f^n(b) <_ {\mathbf E} f^n(a)  $, for some
$n $, hence we cannot have 
$b _1 \leq_O b$, since otherwise
$ f^n(b_1) \leq_ {\mathbf E} f^n(b) <_ {\mathbf E} f^n(a)  $,
contradicting 
$a <_O b_1$.
Hence $b  <_O b_1$.
Here we have used the assumption that
$ \leq_ {\mathbf B} $ is a linear order on $B$.

(vi)
If $b <_O b_1 <_O a$, then
$ f^n(b_1) <_ {\mathbf E} f^n(a)  $, for some
$n $, hence 
$ f^n(b) \leq_ {\mathbf E} f^n(b_1) <_ {\mathbf E} f^n(a)  $, thus
$b <_O a$.

The remaining cases (three elements in $A$ 
or three elements in $B$) are trivial.

We have showed that, for each component $O$, 
the relation $<_O$  given by \eqref{alal} - \eqref{alaluu}
  linearly (strict) orders  $O$.
By the considerations at the beginning, in particular, by (*),
 if, for $d, e \in D$,  we let
\begin{equation} \labbel{ord}  
  \text{$d \leq e$  \quad if  \quad either   $ d \leq_ {\mathbf E} e $,
or $d <_O e$, for some component $O$,}
 \end{equation}
then we get a linear order on $ D$
which extends  $ \leq_ {\mathbf E} $,
thus $(D, \leq )$
amalgamates  
$(A, \leq_ {\mathbf A})$
and $(B, \leq_ {\mathbf B})$ over
$(C, \leq_ {\mathbf C})$ in the class of linear orders.
It remains to show that 
$f$ is order preserving on $D$
with respect to $\leq$.  

As in the proof of Corollary \ref{oper}, from \eqref{al}
and from the assumption
that $f$ is order preserving on both
$\mathbf A$ and $\mathbf B$ it follows    
that if $ d \leq_ {\mathbf E} e $,
then $ f(d) \leq_ {\mathbf E} f(e) $,
hence $ f(d) \leq  f(e) $.
This case covers also \eqref{alal}  and \eqref{alalu}, hence
we only need
to consider the
case in \eqref{alaluu}.
Suppose that 
$a \in A $ and $b \in B$  belong to the same component $O$. 
There are two cases with various subcases.

\emph{Case $b < _O a$.} If $b < _O a$, then 
$ f^n(b) <_ {\mathbf E} f^n(a)  $, for some
$n \geq 1$, by \eqref{alaluu}.
(i) If $n=1$, then $ f(b) <_ {\mathbf E} f(a)  $,
hence $ f(b) < f(a)  $ and we are done.
(iia) If $n>1$ and 
$ f(b) $, $ f(a)  $
lie in the same 
component $O'$, then we get 
$ f(b)  < f(a)  $ applying \eqref{alaluu}
to $<_{O'}$  with $n-1$ in place of  $n$. 
(iib) If $n>1$ and $ f(b) $ and  $ f(a)  $
lie in distinct 
components, then 
either 
$ f(b) \leq _ {\mathbf E} c  \leq _ {\mathbf E} f(a)  $,
or
$ f(a) \leq _ {\mathbf E} c  \leq _ {\mathbf E} f(b)  $,
for some $c \in C$.
But the latter eventuality cannot occur, 
since it implies
$ f^n(a) \leq _ {\mathbf E} f^{n-1}(c)  \leq _ {\mathbf E} f^n(b)  $,
contradicting $ f^n(b) <_ {\mathbf E} f^n(a)  $.
Hence $f(b) \leq f(a)$
in this case, as well.

\emph{Case $a < _O b$.}  
If $a < _O b$, then  
 for no $n$  
$ f^n(b) <_ {\mathbf E} f^n(a)  $, 
hence for no $m$ 
$ f^m(f(b)) <_ {\mathbf E} f^m(f(a))  $.
(i) If  $ f(b) $ and  $ f(a)  $
lie in the same 
component, we are done, by applying  \eqref{alaluu}
with $ f(a) $ and  $ f(b)  $
in place of $a$ and  $b$. 
(ii) If  $ f(b) $ and  $ f(a)  $
lie in distinct components,
then, as above, 
either 
$ f(b) \leq _ {\mathbf E} f(a)  $,
or
$ f(a)  \leq _ {\mathbf E} f(b)  $.
Now observe that if $a < _O b$, then
$ f(b) <_ {\mathbf E} f(a)  $
does not occur, otherwise the first clause in 
 \eqref{alaluu} should have been applied.
Thus $ f(a)  \leq _ {\mathbf E} f(b)  $.

We have proved that $f$ is  
$\leq$-preserving, hence
$\mathbf D=(D, \leq , f)$
is linearly ordered and
amalgamates $\mathbf A$ 
and $\mathbf B$ over  $\mathbf C$.

(b)  We first show that SAP fails. Let $\mathbf C$  be $\mathbb N$
with the usual order and with $f$ interpreted as the  
successor function. 
Let $A= \{ a \} \cup \mathbb N $
and let $\mathbf A$ extend
$\mathbf C$ by setting 
$a < 0$ and $f(a)=0$.
Similarly, let 
 $B= \{ b \} \cup \mathbb N $
with $a \neq b$ 
and let $\mathbf B$ extend
$\mathbf C$ by setting 
$b < 0$ and $f(b)=0$.
If an amalgamating algebra 
  is a linear order and $f$ is still to be
strict order preserving, then $a$ and  $b$ 
  should be  identified, since 
$f(a)=f(b)$ and then 
both $a<b$ and $b<a$
contradict the assumption that 
$f$ is strict order preserving.  
Hence SAP fails.

In order to prove AP,
we shall show that the situation in the above counterexample
essentially provides all kinds of failures of SAP.
In summary, if $a \in A$, $b \in B$ and
$f^n(a)=f^n(b) \in C$, then $a$
and $b$ should be identified.   
After the identification is made,
we are left with a triple
which can be amalgamated
using the techniques of 
part (a) and then it is quite easy to see that
$f$ is strict order preserving
in the amalgamating algebra. 

We now proceed wih the details. Suppose that 
$ \iota _{\mathbf C, \mathbf A} \colon \mathbf C \rightarrowtail \mathbf A$
and 
 $ \iota _{\mathbf C, \mathbf B} \colon \mathbf C \rightarrowtail \mathbf B$
are embeddings.
Let 
\begin{align*} \labbel{ca} 
C_A = \{  \, & a \in A \mid \text{there are $n \in \mathbb N$,
$c \in C$ and  $b \in B$ such that} 
\\
& \text{$f_{\mathbf A}^n(a)=\iota _{\mathbf C, \mathbf A}(c)$
 and $f_{\mathbf B}^n(b)= \iota _{\mathbf C, \mathbf B} (c)$} \, \} 
\text{ and, symmetrically,}  
\\
C_B = \{  \, & b \in B \mid \text{there are $n \in \mathbb N$,
$c \in C$ and  $a \in A$  such that} 
\\
& \text{$f_{\mathbf A}^n(a)=\iota _{\mathbf C, \mathbf A}(c)$
 and $f_{\mathbf B}^n(b)= \iota _{\mathbf C, \mathbf B} (c)$} \, \}. 
   \end{align*} 
 
Notice that $C_A$
and  $C_B$
are closed under applications of $f$,
hence they are domains for  substructures 
 $\mathbf C_A$ and $\mathbf C_B$  of 
 $\mathbf A$  and $\mathbf B$, respectively.  
Moreover,
  $n=0$ is allowed,  
hence $\mathbf C$ embeds both 
in  $\mathbf C_A$ and $\mathbf C_B$.
Now observe that every strict order preserving
unary operation on some linearly ordered set is injective.
Hence if $a \in C_A$,
then the $b \in B$ witnessing that
$a$ satisfies the defining condition for $C_A$ 
is unique. Moreover, such  $b$
clearly belongs to $ C_B$, 
as witnessed, in turn, by $a$.
Thus if we define 
$\varphi \colon C_A \to C_B$
by setting $\varphi(a) =b$,
for $a,b$ as above, we get a 
bijective correspondence 
from $C_A$ onto $C_B$.   
Surjectivity of $\varphi$  is given
   by the symmetrical 
   argument.

The correspondence $\varphi$  
is a homomorphism with respect to $f$,
since if $\varphi(a)=b$, then
$\varphi(f(a))= f(b)$.
This is obvious
if $a \in C$;
otherwise
consider $n-1$
in place of $n$
in the definitions of   $C_A$ and $C_B$.
We are going to show that 
$\varphi$  is also an order isomorphism,
hence 
$\mathbf C_A$ and $\mathbf C_B$
are isomorphic structures.
Indeed, suppose that 
$a, a_1 \in C_A$ with
$f_{\mathbf A}^n(a)=\iota _{\mathbf C, \mathbf A}(c)$
and 
$f_{\mathbf A}^m(a_1)=\iota _{\mathbf C, \mathbf A}(c_1)$
with, say, $m>n$.
If $r=m-n$,
then   
$f_{\mathbf A}^m(a)=\iota _{\mathbf C, \mathbf A}(f^r_{\mathbf C}(c))$
and, since $\iota _{\mathbf C, \mathbf A}$
is an embedding and $f$ is strict order preserving, we get
that $a < a_1$ if and only if 
$f_{\mathbf A}^m(a) < f_{\mathbf A}^m(a_1)$,
if and only if    
$f^r_{\mathbf C}(c) < c_1$. 
The last inequality is computed in 
$\mathbf C$, hence it is also equivalent to 
$\varphi (a) < \varphi(a_1)$.

Since $\mathbf C_A$ and $\mathbf C_B$
are isomorphic,
then, by 
replacing 
$\mathbf A$, $\mathbf B$  and $\mathbf C$
with suitable isomorphic structures,
it is no loss of generality to assume that
$\mathbf C_1 \subseteq \mathbf A, \mathbf B$,
where $\mathbf C_1$
is isomorphic to 
$\mathbf C_A$. Since 
$\mathbf C$ embeds into   $\mathbf C_1$,
if we can amalgamate
the copies of  $\mathbf A$ and
$\mathbf B$ over  
$\mathbf C_1$, then we have embeddings amalgamating
the original  structures.
The relevant property of 
 $\mathbf C_1$
that we have obtained is that
\begin{equation}\labbel{*}    
 \text{if $a \in A$, $b \in B$,
 $c \in C_1$ and $f_{\mathbf A}(a)=c=f_{\mathbf B}(b)$,
 then $a=b \in C_1$.}  
  \end{equation}

Now apply the construction in (a)
with $\mathbf C_1$ in place of $\mathbf C$.
Since $f$ is strict order preserving  on $\mathbf A$
and $\mathbf  B$, it is in particular
order preserving, hence (a) can be applied, obtaining some
structure $\mathbf D$ with an order preserving $f$.  
It remains to show that $f$ is strict order preserving, and we shall show
that this follows from \eqref{*}. 
So let $d<e$. We know from (a)
that $f(d) \leq  f(e)$; it remains to show that
$f(d) \neq  f(e)$.
If either $d,e \in A$
or  $d,e \in B$,
this is immediate from the assumption
that $f$ is strict order preserving on $\mathbf A$,
 respectively, $\mathbf  B$.
Otherwise, say, $d \in A$ and $e \in B$,
thus   $f(d) \in A$ and $f(e) \in B$.
If $f(d)=f(e)$, then  
$f(d), f(e) \in C_1 = A \cap B$,
hence $d=e$, by \eqref{*}.  
This contradicts $d < e$.

(c) 
We shall present two counterexamples, since they 
have quite distinct features.

(c)(i) Let $C= \{  0 \} $ with the only possible interpretations,
$A= \{ a, 0\}$, $B= \{ b, 0\}$,
$A\cap B=C$ with 
\begin{equation}\labbel{c}
\begin{aligned} 
&\text{$ a < _{\mathbf A} 0$,
$ f _{\mathbf A}(a) =  0$,
$ h _{\mathbf A}(a) =  a$, and}
\\ 
&\text{$ b <_{\mathbf B} 0$,
$ f _{\mathbf B}(b) =  b$,
$ h _{\mathbf A}(b) =  0$.}
\end{aligned}   
 \end{equation}    

If $\mathbf A$, $\mathbf B$ and $\mathbf C$ 
 can be amalgamated to a \emph{linear} order, then in 
the amalgamating algebra we have either 
$a \leq b$ or $b \leq a$.
If $f$ is required to be
order preserving, the first eventuality cannot occur, since then 
$0 = f(a) \leq f(b) = b$ and similarly the second eventuality cannot 
occur, if $h$  is required to be
order preserving.

(c)(ii) If we want $f$ and $h$ to be strict
order preserving, consider 
 $\mathbb N$, as $\mathbf C$, with the standard order and interpret both
$f$ and $h$ in $\mathbf C$ as the successor function.
Then extend $\mathbf C$
with one more element in two possible ways 
with $A= \{ a \} \cup  \mathbb N $, 
$B= \{ b \} \cup  \mathbb N $
and the same new relations \eqref{c} as above.
Then repeat the same argument.  

(d)  In this case the order used in the proof 
of Corollary  \ref{losp} works.
In detail, given  
 $\mathbf C \subseteq \mathbf A, \mathbf B$
to be amalgamated, let 
$D= A \cup B$ and, for each $f \in  F$,  define $f$ on $D$ 
in the unique compatible way.
Thus each $f$ is bijective, since $f$  is bijective both on $A$
and $B$ and the values of $f$ agree on 
$C= A \cap B$.   
Then extend
the order $\leq_{ \mathbf E}$ on $D$  
from Theorem \ref{FJ} by setting 
$d \leq e$  if 
\begin{equation}\labbel{eqeq}        
\begin{aligned}
&\text{either } 
d \leq_{ \mathbf E} e, \text{ or }   
\\
&\text{not $d \leq_{ \mathbf E} e$, 
not $e \leq_{ \mathbf E} d$ and $d \in A$,  $e \in B$.} \quad \qquad
 \end{aligned}
 \end{equation}

By the arguments in (a), if the second alternative in equation \eqref{eqeq} holds,
then $d$ and $e$ belong to the same component.
Now notice that if $f$ is bijective and order preserving
on $C$ and 
$(C_1, C_2)$ is a cut of $\mathbf C$, then
$(f(C_1), f(C_2))$  is a cut.
The assumption that $f$ is surjective is used
in order to get  $f(C_1) \cup  f(C_2)= C$.
Hence  if $d$ and $e$ belong to the  component 
$O$ associated to $(C_1, C_2)$, 
then   $f(d)$ and $f(e)$ belong to the component
$O'$ associated to $(f(C_1), f(C_2))$, since
$f$ is strict order preserving on both $\mathbf A$ and  $\mathbf B$. 

Now we can prove that each $f$
is order preserving with respect to $\leq$. 
If  
$d <_{ \mathbf E} e$,
then we get 
$f(d) \leq_{ \mathbf E} f(e)$.
On the other hand, if $d < e$
is given by the second alternative 
in \eqref{eqeq}, then $d \in A$,
$e \in B$ and $d,e $ belong to the same component,
call it $O$.
But also   $f(d) \in A$,
$f(e) \in B$ and $f(d),f(e) $ belong to the same component 
$O'$ as described above,
hence \eqref{eqeq} gives 
$f(d) < f(e)$.
Put in another way, the first alternative in
\eqref{alaluu} never occurs when dealing
with bijective functions. 

The conclusion follows from the fact that a bijective order
preserving function on a linearly ordered set is necessarily an 
order automorphism.
\end{proof} 

In the special case of a single strictly increasing automorphism
AP in Theorem \ref{los} is a consequence of \cite[Theorem 2.2]{LP}.
Recall that a theory with model completion has
AP.  

\begin{remarks} \labbel{ultt}
The  counterexample in (c)(i)
in the proof of Theorem \ref{los}  shows that
 the class of \emph{finite} linearly ordered sets
with two order preserving operations
fails to have AP. Indeed, the counterexample shows a bit more.

Recall that 
 an \emph{(order-theoretical) closure operation} on some poset $P$ 
is an order preserving unary operation  $f$ such that 
  $f(f(x))=f(x) \geq x$ holds for every $x \in P$.
See \cite{Er} for information about closure operations, pictures
and for the interest of the notion in the general
order-theoretical setting. 
The  counterexample in (c)(i)
  shows that
 the class of (finite) linearly ordered sets
with two closure  operations
fails to have AP.

The  counterexample (c)(ii) works both for the strict 
and the nonstrict case, 
but in the former situation the counterexample should necessarily be
infinite, since a  strict order preserving operation
 is the identity on a finite 
linearly ordered set.

Actually, the example shows that it is not always the case that
a triple of linearly ordered set with 
two strict order preserving 
operations can be amalgamated into a
linearly ordered set with 
two  order preserving 
operations, namely, without requiring 
in the amalgamating structure  that the operations
are strict order preserving.

The proof of Theorem \ref{los}(d)
shows that if $\mathbf C$
is a linearly ordered set with
an $F$-indexed set of 
automorphisms,
then $\mathbf C$ is a strong amalgamation base
for the class of   linearly ordered sets with
an $F$-indexed set of 
strict order preserving unary operations.
 \end{remarks}

\section{Order reversing operations} \labbel{r}  

If we consider order reversing operations, 
the arguments of the previous section generally carry over.
However, a linearly ordered set with an order reversing
operation has at most one element
$c$ such that $g(c)=c$, and obviously embeddings
must preserve such ``centers'', if they exist.
This fact prevents strong amalgamation.
Moreover, elements greater than the center
should be treated in a different---but symmetrical---way
 in comparison with elements smaller than the center.
 We now give precise definitions and collect some
trivial facts about these notions.

\begin{definition} \labbel{cent}
If $\mathbf C$ is a linearly ordered set 
with a unary operation
$g$, an element $c$ of $ C$ is said to be a \emph{center},
or a \emph{fixed point} 
of  $g$ if $g(c)=c$.

An element $d$ of $ C$
is an \emph{upper} (resp., \emph{lower}) 
element if $g(d) <d$ (resp., $g(d)>d$).
 \end{definition}

\begin{lemma} \labbel{lem}
Suppose that $\mathbf C$ is a linearly ordered set 
with one order reversing unary operation
$g$.
  \begin{enumerate}[(a)]    
\item
Every element of $C$ is either upper,  lower,
or a center. The alternatives are mutually exclusive 
and $\mathbf C$ has at most one center.
\item
Suppose that   $\mathbf C$ has a center $c$.
Then,  for every $d \in C$,
$d$ is upper if and only if $ c < d$
and  $d$ is lower if and only if $ d <c $.
\item
If $d$ is upper, then $g(d)$ is either lower or the
center, and symmetrically.   
All the upper elements are greater
than all the lower elements,
\item
If   $\mathbf C$ has no center,
then  $\mathbf C$ 
can be extended by adding just one element 
(in a unique way
modulo isomorphisms preserving $\mathbf C$)
 to a 
linearly ordered set $\mathbf C^*$
with an order reversing unary operation
and a center.
If $g$ is strict order reversing in $\mathbf C$,
then the operation in $\mathbf C^*$ is strict order reversing, too. 

\item
Suppose that 
$\iota \colon \mathbf C \rightarrowtail \mathbf A$
is an embedding.

If $\mathbf C $ has a center $c$,
then $  \mathbf A$ has a center and
$\iota (c)$ is the center of $  \mathbf A$.

 If $\mathbf C $ has not a center and 
$  \mathbf A$ has a center, then
$\iota$ extends uniquely to an embedding 
from 
$\mathbf C^*$ to 
$  \mathbf A$, 
where $\mathbf C^*$
is defined as in (d).

 If neither  $\mathbf C $ nor 
$  \mathbf A$ have a center, then
$\iota$ extends uniquely  to an embedding 
from 
$\mathbf C^*$ to 
$  \mathbf A^*$.
\end{enumerate} 
 \end{lemma}

 \begin{proof} 
(a) - (c) are immediate from the assumptions
that  $\mathbf C$ is linearly ordered and
$g$ is order reversing. For example, 
to prove (c), observe that if $g(d) <d$, then 
$g(g(d)) \geq g(d)$, hence either 
$g(g(d)) > g(d)$, thus $g(d)$ is lower, or 
 $g(g(d)) = g(d)$, thus $g(d)$ is the center.
Suppose that $d$ is upper and  $e$ is lower. If $d < e$, 
then $g(d) <d< e < g(e)$
contradicts the assumption that 
$g$ is order reversing. Hence $e< d$, since the order is linear
and $d$ and $e$  are necessarily distinct. 

(d) By (b), the new element supposed to be a center, call it  $c$,
should be greater than all the lower elements and smaller
then all the upper elements, hence the position of 
$c$ in the order is fully determined.
Setting $g(c)=c$, clause 
(c) implies that, endowed with the  above structure,
$C \cup \{  c\} $ is     
linearly ordered and $g$ 
is order reversing.

(e) It follows from the definitions that an embedding
(actually, just a morphism) sends a center to a center.
All the rest follows from (a) - (d).
\end{proof}

It follows from Lemma \ref{lem}
that embeddings preserve 
upper 
and lower elements, as well as
centers,
if they exist. In particular,
given a triple to be amalgamated,
there is no need to mention some specific structure
$\mathbf A$, $\mathbf  B$ or $\mathbf  C$, when
referring to the center.

\begin{thm} \labbel{rev}
(a) The classes  
$\mathcal {LO}_{r}$, resp.,
$\mathcal {LO}_{sr}$
 of linearly ordered sets 
with one order reversing, resp., 
one strict order reversing unary operation
 have  AP but not SAP.

(b)
The classes of linearly ordered sets with 
two order reversing, resp.,
two strict order reversing unary operations
have not AP.
Similarly for the case of an order preserving and
an order reversing operation.
AP fails even if we assume that all the operations have a 
common center. 
\end{thm}

\begin{proof}
(a) To prove that SAP fails, just let 
 $\mathbf C$ have no center and 
 $\mathbf A$, $\mathbf B$ have a center.
In any amalgamating structure the centers of 
$\mathbf A$ and $\mathbf B$ must be identified,
by Lemma \ref{lem}, hence SAP fails.
The simplest concrete example is when
$\mathbf C$ is an empty model and 
$\mathbf A$ and $\mathbf B$
have only one element, necessarily, the center.
However, we want to prove that
also the weaker version of SAP fails 
when $\mathbf C$ is required to be nonempty.
Cf.\ the final comment in Definition \ref{sap}. 

So let $\mathbf C $ be the model with domain $C= \{  - \infty, \infty\} $
and such that  $ - \infty < \infty$, 
$g (- \infty) = \infty$
and  $g ( \infty) = - \infty$.
Extend $\mathbf C$ to $\mathbf A$ 
by letting $ A= \{  - \infty, a,  \infty\} $,
with $ - \infty < a < \infty$ 
 and $g (a) = a$. Similarly,
extend $\mathbf C$ to $\mathbf B$ 
by letting $ B= \{  - \infty, b,  \infty\} $
with $b \neq a$, $ - \infty < b < \infty$ 
 and $g (b) = b$
(in fact, $\mathbf A$ and $\mathbf B$ 
are just two isomorphic copies of the structure 
$\mathbf C^*$ constructed in Lemma \ref{lem}(d)).
By Lemma \ref{lem},
in any amalgamating structure we must have
$a=b$, hence SAP fails.
Notice that $g$ is strict order reversing on 
$\mathbf A$, $\mathbf B$ and $\mathbf C$,
but $a=b$ in any amalgamating structure $\mathbf D$,
even if $g$ is assumed to be (possibly, not necessarily strict) order
reversing in $\mathbf D$.   

The proof of AP is similar to the proof 
of Theorem \ref{los}(a)(b), except that
the possible
 overlapping of 
centers should be fixed (in the case of
order reversing operations this is the only obstacle
to strong amalgamation) and that the actual definition of the linear
order involves still another division into cases. 

 Let $\mathbf C \subseteq \mathbf A, \mathbf B$
be a triple to be amalgamated.
If some  algebra above has
no center, add a center to it
according to Lemma \ref{lem}(d). 
Possibly, replace $\mathbf A$ and $ \mathbf B$ 
with isomorphic copies, so that their centers
are identified (this is necessary exactly in case the 
original $\mathbf A$ and $ \mathbf B$ have some center
and $\mathbf C$ has  not a center).
 Because of Lemma \ref{lem}, there is just one way to add
the centers and the original embeddings
can be extended in a unique way.

Hence we can suppose that 
$\mathbf A$, $\mathbf B$ and $ \mathbf C$
all have a center.
As in the proof of  Theorem \ref{los},  Theorem \ref{FJ}
and Corollary \ref{oper}  
furnish an amalgamating  \emph{partial} order
$\mathbf E=(D, \leq _ {\mathbf E} , g)$  
with an order reversing  operation.
Lengv\'{a}rszky \cite{L}
showed that a poset with a unary order reversing function
$g$ can be linearized in such a way that $g$
is still order reversing if and only if    
$g^2$ is acyclic and $g$  has at most one fixed point
(in the original poset).  Hence we can apply
\cite{L} in order to get a proof of the positive part of (a).
As in the case of Theorem \ref{los},
we shall present a direct and more explicit construction.  

To simplify the notation, let 
$<_ {\mathbf E}^n$
be $<_ {\mathbf E}$
for $n$ even, and    
$<_ {\mathbf E}^n$
be $>_ {\mathbf E}$
for $n$ odd.
Apply a similar convention for $\leq^n _ {\mathbf E}$.
Recall from the proof of Theorem \ref{los} 
the definition of a component, and recall that
the relative order between two elements lying
in distinct components is completely determined by
$\leq _ {\mathbf E}$. Notice also that, as we mentioned,
the definition of the components does not rely
on the operations, it depends only on the orderings.

As in the proof of \ref{los}(a),
we need to linearize each component.
 Let us call a component \emph{lower}
if its elements are  $\leq _ {\mathbf E} c$
and \emph{upper} otherwise, where $c$
is the center of $\mathbf C$.
The distinction makes sense, since each
component is convex and contained in 
$  D \setminus C$.    
If $O$ is a lower component, linearize 
$O$ according to the conditions \eqref{alal} - \eqref{alaluu},
but replacing    
$ f^n(b) <_ {\mathbf E} f^n(a)  $ in \eqref{alaluu}
by $ g^n(b) <^n_ {\mathbf E} g^n(a)  $. 
The proof that $<_O$ is a linear order on $O$
carries over just by replacing   
  $<_ {\mathbf E}$
and   $\leq_ {\mathbf E}$, respectively, 
by
$<_ {\mathbf E}^n$
and
$\leq_ {\mathbf E}^n$
in all the expressions involving $f^n$
(here, $g^n$)
and with no further modification.

As far as upper components are concerned,
we have to exchange the role of $A$ and $B$ 
in \eqref{alaluu}, since we want $g$ to be order reversing.  
In detail, replace \eqref{alaluu} by 
\begin{gather}  
\labbel{alaluuu} 
\begin{aligned}
&\left. \begin{aligned}  
 &a <_O b, 
&&\text{if $ g^n(a) <^n_ {\mathbf E} g^n(b)  $},
\\[- 2pt]&&&\text{for some
$n \geq 1$, or} 
\\
& b<_O a, &&\text{otherwise,}
\end{aligned} 
\ \ \right\}
&&\text{if $a \in A \cap O$ and $b \in B \cap O$}
\\
\end{aligned}
\\
\labbel{bup}
\tag*{}
 \text{  for upper components } 
\end{gather} 
By symmetry, $<_O$ is a linear order
on $O$ in this case, too.

We now can define $\leq$ on $D$ 
as in equation \eqref{ord}.
The proof that $g$ is order reversing
with respect to $\leq$ is similar to 
\ref{los}(a), just considering separately
the cases when   $a$ and  $b$
belong to the same lower or upper 
component.
For example, we shall treat the case when
$a \in A $ and $b \in B$  belong to the same upper component $O$. 

\emph{Case $a < _O b$.} 
If $a < _O b$, then 
$ g^n(a) <^n_ {\mathbf E} g^n(b)  $, for some
$n \geq 1$, by \eqref{alaluuu}.
(i) If $n=1$, then
$ g(a) <^1_ {\mathbf E} g(b)  $,
namely,
$ g(b) <_ {\mathbf E} g(a)  $,
by the definition of 
$<^n_ {\mathbf E}$, 
hence $ g(b) < g(a)  $ and we are done.
(ii) If $n>1$, 
then
$ g^{n-1}(g(a)) <^n_ {\mathbf E} g^{n-1}(g(b))  $,
that is, 
$ g^{n-1}(g(b)) <^{n-1}_ {\mathbf E} g^{n-1}(g(a))  $.
(iia) First suppose that $ g(b) $ and $ g(a)  $
lie in the same 
component $P$.
Then $P$ is a lower component, by 
Lemma \ref{lem}(c)
applied to $\mathbf A$ and $\mathbf B$.
Notice that, by construction, each component 
has empty intersection with $C$, hence 
$ g(b) $ and $ g(a)  $ are not  centers,
since $C$ is assumed to have a center and the center is unique. 
Thus we get 
$ g(b)  < _P g(a)  $ applying the modified version of \eqref{alaluu}
with $n-1$ in place of  $n$. 
(iib) If $ g(b) $ and  $ g(a)  $
lie in distinct 
components, then 
either 
$ g(b) \leq _ {\mathbf E} c  \leq _ {\mathbf E} g(a)  $,
or
$ g(a) \leq _ {\mathbf E} c  \leq _ {\mathbf E} g(b)  $,
for some $c \in C$.
The latter eventuality cannot occur, 
since it implies
$ g^n(a) \leq^{n-1} _ {\mathbf E} g^{n-1}(c)  \leq^{n-1} _ {\mathbf E} g^{n}(b)  $,
hence
$ g^n(b) \leq^{n} _ {\mathbf E}  g^{n}(a)  $
contradicting $ g^n(a) <^n_ {\mathbf E} g^n(b)  $.
Hence $g(b) \leq g(a)$
in this case, as well.

\emph{Case $b < _O a$.} 
If $b < _O a$, then  
 for no $n$  
$ g^n(a) <^n_ {\mathbf E} g^n(b)  $, 
hence for no $m$ 
$ g^m(g(a)) <^{m+1}_ {\mathbf E} g^m(g(b))  $,
that is, 
for no $m$ 
$ g^m(g(b)) <^{m}_ {\mathbf E} g^m(g(a))  $.
(i) If  $ g(b) $ and  $ g(a)  $
lie in the same 
component $P$, necessarily, as we mentioned, a lower component,
then we get $g(a)<_Pg(b)$, by applying 
the variant of \eqref{alaluu}
with $ g(a) $ and  $ g(b)  $
in place of $a$ and  $b$.
Hence  $g(a)<g(b)$.
(ii) If  $ g(b) $ and  $ g(a)  $
lie in distinct components,
then, as above, 
either 
$ g(b) \leq _ {\mathbf E} g(a)  $,
or
$ g(a)  \leq _ {\mathbf E} g(b)  $.
If  $b < _O a$, then
$ g(a) <^1_ {\mathbf E} g(b)  $,
equivalently,
$ g(b) <_ {\mathbf E} g(a)  $,
does not occur, otherwise the first clause in 
 \eqref{alaluuu} should have been applied.
Thus $ g(a) \leq _ {\mathbf E} g(b)  $,
hence $ g(a) \leq  g(b)  $.

The case when $g$ is assumed to be
strict order reversing presents no
essential difference with respect to the proof of
Theorem \ref{los}(b).
The only minor detail is in the proof that $\varphi$,
as defined in the proof of \ref{los}(b),
is an order isomorphism.
In the present case, assume that
$a, a_1 \in C_A$,
$g_{\mathbf A}^n(a)=\iota _{\mathbf C, \mathbf A}(c)$
and 
$g_{\mathbf A}^m(a_1)=\iota _{\mathbf C, \mathbf A}(c_1)$
with, say, $m>n$.
If $r=m-n$,
then   
$g_{\mathbf A}^m(a)=\iota _{\mathbf C, \mathbf A}(g^r_{\mathbf C}(c))$.
Since $g$ is strict order reversing, then
$a < a_1$ if and only if 
$g_{\mathbf A}^m(a) <^m g_{\mathbf A}^m(a_1)$,
if and only if    
$g^r_{\mathbf C}(c) <^{m} c_1$. 
As in \ref{los}(b), the last inequality is computed in 
$\mathbf C$, hence it is also equivalent to 
$\varphi (a) < \varphi(a_1)$.
All the rest goes as in \ref{los}(b)

(b)(i) We provide the example of three
finite nonamalgamable algebras with two 
order reversing
operations $g$ and $k$.
Let $C= \{  c \} $;
$A= \{ a, c, d\}$,
with $ a < _{\mathbf A} c < _{\mathbf A}  d$,
$ g _{\mathbf A}(a) = g _{\mathbf A}(d) =  c$,
$ k _{\mathbf A}(a) =  d$,
$ k _{\mathbf A}(d) =  a$,
 and 
$B= \{ b, c, e\}$,
with $A \cap B = \{ c \} $, $ b< _{\mathbf B} c < _{\mathbf B}  e$,
$ g _{\mathbf B}(b) =  e$,
$ g _{\mathbf B}(e) =  b$,
$ k _{\mathbf B}(b) = k _{\mathbf B}(e) =  c$,
Then argue as in \ref{los}(c)(i). 

(b)(ii)
In this example we construct
  three
nonamalgamable finite algebras with  an
order preserving operation
$f$ and an 
order reversing
operations $g$.
Let the domains $A$, $B$, $C$, the orderings
 and the operations $g$ be as in (b)(i).
Let  $f_{\mathbf A} $ be the identity  and
$f_{\mathbf B} (b) = f_{\mathbf B} (e) = c$.
If $a  \leq b$ in some amalgamating structure
with  $g$ order reversing,    
then  $e= g(b) \leq g(a) =c$, contradicting $c < e$.   
If $b \leq a$ in some amalgamating structure with
$f$ order preserving, then
$c = f(b) \leq f(a) = a$, again a contradiction.

(b)(iii) Now we present  three
nonamalgamable algebras with two 
strict order reversing
bijective operations $g$ and $k$
and which cannot be amalgamated
into a linear order on which
$g$ is order reversing. 
Let $C= \{  - \infty, \infty\} $
with  $ - \infty < \infty$, 
$g_{\mathbf C} (- \infty) =k_{\mathbf C} (- \infty) = \infty$
and  $g_{\mathbf C} (\infty) =k_{\mathbf C} ( \infty) =- \infty$.
Extend $\mathbf C$ to $\mathbf A$ 
by letting $ A= C \cup \mathbb Z$,
with $ - \infty < z < \infty$
and $g_{\mathbf A} (z) =k_{\mathbf A} (z) = -z$,
 for every $z \in \mathbb Z$.
Let $\mathbb Z'= \{  \dots, -2', -1',0', 1', 2', \dots  \} $
be a disjoint copy of $\mathbb Z$ and 
 extend $\mathbf C$ to $\mathbf B$ 
by letting $ B= C \cup \mathbb Z'$,
with $ - \infty < z' < \infty$,
 $g_{\mathbf B} (z') = -z'$ and
$k_{\mathbf B} (z')=  -z' +2'$,
 for every $z' \in \mathbb Z'$.

In view of Lemma \ref{lem},
in any amalgamating structure
with $g$ (not necessarily strict) order reversing,
the centers
 $0$ and $0'$ with respect to $g$
  should be identified,
but this is incompatible with
 $k_{\mathbf A} (0)=  0$
and $k_{\mathbf B} (0')=  2'$.
 
Notice that $k_{\mathbf A} (0)=  0$
and $k_{\mathbf B} (0') \neq 0'$ are the only 
properties of $k$ needed in  the above 
argument, hence, by changing the other
values of $k$,  the counterexample can be modified 
in order to take care of the case of a (strict) order reversing 
together with a
(strict) order preserving operation,
possibly both bijective. Actually, there are 
plenty of further similar possibilities.
 
(b)(iv) The main point in (b)(iii)
above is that the operations in $\mathbf  C$ have no center
and then centers are added in different ways to $\mathbf A$ and
$\mathbf  B$. On the other hand, we can merge the ideas in 
(b)(i) and \ref{los}(c)(ii) in order to get failure of AP
even under the assumption that  the two operations
are strict order reversing with a common center. 

So let $\mathbf  C$ be $\mathbb Z \setminus \{ -1, 1 \} $
with the standard order, $g_{ \mathbf  C} (0) = k_{ \mathbf  C} (0)=0$,
$g_{ \mathbf  C} (n) = k_{ \mathbf  C} (n)=-n$
and 
$g_{ \mathbf  C} (-n) = k_{ \mathbf  C} (-n)=n+1$,
for $n \in \mathbb N \setminus \{ 0,1 \}  $.
Extend $\mathbf  C$ to $\mathbf A$ 
with $A= \mathbb Z$,  
 $g_{ \mathbf  A} (1) =   k_{ \mathbf  A} (1)=-1$
and 
$g_{ \mathbf  A} (-1) = 1$, $  k_{ \mathbf  A} (-1)=2$.
Let $B=C \cup \{ 1', -1' \} $ 
with $-2<-1'<0<1'<2$, 
$g_{ \mathbf  B} (1') =   k_{ \mathbf  B} (1')=-1'$
and 
$g_{ \mathbf  B} (-1') = 2$, $  k_{ \mathbf  B} (-1')=1'$.
 In any amalgamating structure with a linear order, either
$1 \leq 1'$ or $1' \leq 1$.
If $1 \leq 1'$, then
$2=k^2(1) \leq k^2(1')=1' $, a contradiction. 
If $1' \leq 1$, then
$2=g^2(1') \leq g^2(1)=1 $, still a contradiction. 

If we want a counterexample with a 
strict order preserving operation $f$
and a strict order reversing operation $g$,
again, with a common center,
just take $f=k^2$ in the above counterexample.   
\end{proof}

\begin{remark} \labbel{tutdue}
The proof of Theorem \ref{rev}(a) 
shows that 
the class of linearly ordered sets with 
an order reversing unary operation with a center
has SAP.

Actually, in the class of linearly ordered sets with 
an order reversing unary operation,
a structure  $\mathbf C$ is a
strong amalgamation base if and only if 
$\mathbf C$ has a center.
 \end{remark}

The counterexample (b)(iii)
in the proof of Theorem \ref{rev}
shows that Theorem \ref{los}(d), as it stands,
does not generalize to order reversing bijective operations,
\emph{antiautomorphisms}, for short. 

The counterexamples (b)(i) and (b)(iv)
show that the class of linearly ordered sets
with two order reversing operations with the same center
fails to have  AP.

However, Theorem \ref{los}(d) does generalize if we put together
the two assumptions. 
Moreover, we can deal with  automorphisms
and antiautomorphisms at the same time, provided they all respect the same center.

\begin{theorem} \labbel{td}
For every pair $F$ and $G$ of sets,
let $\mathcal {LO}_{FGac}$ be the class of  linear
orders with an $F$-indexed 
family of automorphisms and 
a $G$-indexed 
family of antiautomorphisms 
such that all the operations in $F$ and in $G$ 
have a common center.
Then  $\mathcal {LO}_{FGac}$ has SAP.
 \end{theorem}

\begin{proof}
If $G= \emptyset $, this
is Theorem \ref{los}(d); actually, no assumption on centers is needed.
So let us assume that  $G \neq \emptyset $, hence  the center is unique and
is preserved by embeddings, by Lemma \ref{lem}.
  
Given a triple $\mathbf A$, $\mathbf B$, $\mathbf C$
to be amalgamated and with center $c$,  extend
the order $\leq_{ \mathbf E}$    
given by Theorem \ref{FJ} on $A \cup B$ by setting 
$d \leq e$  if  either
\begin{align} \labbel{eq1}  
&
d \leq_{ \mathbf E} e, \text{ or }   
\\ \labbel{eq2} &\text{not $d \leq_{ \mathbf E} e$, 
not $e \leq_{ \mathbf E} d$, and $d \in A$,  $e \in B$, $d,e <_{ \mathbf E} c$, or}
\\ \labbel{eq3}
 &\text{not $d \leq_{ \mathbf E} e$, 
not $e \leq_{ \mathbf E} d$, and $d \in B$,  $e \in A$, $c <_{ \mathbf E} d,e$.}  
\end{align}
The definition provides a linear order since if, say,
$d <_{ \mathbf E} c <_{ \mathbf E} e$, then
$d <_{ \mathbf E} e$, hence \eqref{eq1} applies.
If $f_{ \mathbf A}$ 
and $f_{ \mathbf B}$ are 
automorphisms of $ \mathbf A$ 
and $\mathbf B$,
define $f$ on $A \cup B$  
in the unique compatible way as in Corollary \ref{oper}. 
If   $d \leq_{ \mathbf E} e$,
then $f(d) \leq_{ \mathbf E} f(e)$ as in \ref{oper}.
Since $f_{ \mathbf A}$ 
and $f_{ \mathbf B}$ have 
 center $c \in C$,
then   $d,e <_{ \mathbf E} c$
implies  $f(d), f(e) <_{ \mathbf E} f(c) =c$,
hence the arguments in the proof of Theorem \ref{los}(d)
show that if $d < e$ and clause \eqref{eq1}
does not apply, then 
$f(d)<f(e)$.  If $c <_{ \mathbf E} d,e$, 
the symmetrical arguments apply when 
$ A$ 
and $  B$
are exchanged.
Thus $f$ is order preserving on $D$. 
Since $f$ is bijective, it is an automorphism. 

Now suppose that  $g$ is interpreted as an
antiautomorphism 
on $A$ and $B$ and, again,
define $g$ on $A \cup B$ in the unique compatible way. 
 If
  $d,e <_{ \mathbf E} c$,
then 
 $c = g(c) <_{ \mathbf E} g(d), g(e)$.
Recalling the definition of a component from the proof
of Theorem \ref{los}(a), and by the comments before
(*) there, if 
either \eqref{eq2} or \eqref{eq3} applies, then
$d$ and $e$ lie in the same component.   
Arguing in a way similar to 
\ref{los}(d),
the assumption that 
$g$ is bijective implies that   
if $d$ and $e$ lie in the same component,
say, the component associated to 
the cut $ (C_1, C_2)$,
then $g(d)$ and  $g(e)$ lie in the  component 
associated to
$ (g(C_2), g(C_1)$.
Thus if \eqref{eq2} applies to
$d$, $e$, then 
\eqref{eq3} applies to
$g(e)$, $g(d)$ and conversely.
This implies that $g$
 is order reversing on $D$,
hence an antiautomorphism, since $g$
is bijective. 
\end{proof}

\section{Further remarks} \labbel{fur} 

In this section we present a few model-theoretical consequences of
the above results.
It is almost immediate from Theorems \ref{los} and \ref{rev}
that the classes of finite linearly ordered sets  with 
one order preserving, resp., one order reversing unary operation
have a Fra\"\i ss\'e limit.
Moreover, say, if $T$ is the theory
of  linearly ordered sets  with 
one order preserving unary operation $f$ 
satisfying $f^{m+1}(x)= f^m(x)$,
for some fixed $m$,  and $\mathbf M$
is the Fra\"\i ss\'e limit of the class of finite models of $T$,
then $Th(\mathbf M)$ is $ \omega$-categorical and is the 
model completion of $T$.  

With a bit more notation,
we can prove AP, JEP  and the existence of Fra\"\i ss\'e 
limits for many more classes.
The relevant aspect in the following
considerations is that in all the previous constructions
the amalgamating model has been always constructed on the 
set theoretical union of $A$ and $B$.
We shall elaborate further on this aspect in 
\cite{apu}. 

Recall the definitions of the classes 
  $\mathcal {PO}_F$,
 $\mathcal {LO}_{p}$, $\mathcal {LO}_{sp}$,
$\mathcal {LO}_{Fa}$,
$\mathcal {LO}_{r}$, $\mathcal {LO}_{sr}$ and $\mathcal {LO}_{FGac}$
from Proposition
 \ref{oper} and  Theorems  \ref{los}, \ref{rev}, \ref{td}.
Recall the definition of a closure operation
from Example \ref{ultt}.

For every class 
$\mathcal K$ of structures and every set $H$ of appropriate conditions, 
let $ \mathcal K^H$
denote the  subclass of $\mathcal K$  
consisting of those structures in $\mathcal K$ satisfying 
all the conditions in $H$.
We allow $H$ to be the empty set of conditions;
in this case $ \mathcal K^H = \mathcal K$.
In a few cases, for certain combinations of $\mathcal K$
and $H$,  the class 
$ \mathcal K^H$ will turn out to be an empty
class; formally, the results remain true
in this trivial situation.

For each class we have considered in this note, AP and SAP
are preserved by adding various kinds of conditions.
In some cases, the classes we have considered
have the Joint Embedding Property (JEP),
even when AP  fails.

\begin{lemma} \labbel{inpiu} 
The classes 
  $\mathcal {PO}_F^H$,
 $\mathcal {LO}_{p}^H$, 
$\mathcal {LO}_{Fa}^H$
 and $\mathcal {LO}_{FGac}^H$
have SAP and JEP, for any pair of sets 
$F$ and $G$  and for any set  $H$ of conditions
chosen among the following ones.

%%%%CAMBIATO QUI
The ordered set has no maximum (minimum); is finite;
finitely generated; countable; of cardinality $< \lambda $, for 
$\lambda$ an infinite cardinal; 
is well-ordered;  some operation 
$f$ (or some iteration $f ^ \ell$, $\ell \in \mathbb N$) has some  (no) fixed point;
 is surjective; is (strictly) increasing (decreasing);
is a closure operation;
for some 
$m , n \in \mathbb N$ satisfies $f^{m+1}(x)=f^n(x)$ for every
(some) $x$;
some given pair of operations commute.
In general, we can allow any condition
which can be expressed by a universal-existential
first-order sentence such that only one variable is bounded by
the universal quantifier. 

The classes $\mathcal {LO}_{sp}^H$,
$\mathcal {LO}_{r}^H$, $\mathcal {LO}_{sr}^H$
have AP and JEP, for any set  $H$ of conditions
as above.

The class
of linearly ordered sets with any (fixed in advance) number of 
order preserving and strict order preserving 
unary operations
has JEP. Each subclass determined by any set of
conditions as above has JEP .
\end{lemma}

\begin{proof}
If some  property from $H$ holds
in $\mathbf A$, $\mathbf B$ and  $\mathbf C$,   
then the property holds in the amalgamating structure $\mathbf  D$,
 since in each case we have constructed $\mathbf  D$ 
on $A \cup B$. 
Hence (S)AP holds in all the classes under consideration.

All the classes for which we have proved AP
have also JEP,
since, for languages without constants, JEP means
exactly
that the empty structure
is an amalgamation base. 
Formally, the class
$\mathcal {LO}_{FGac}$
has not an empty model; 
however, modulo isomorphism, there is a unique 
``initial'' $1$-element model,
hence JEP follows from AP.
Notice that it is not necessary to assume that the 
$1$-element model, call it  $\mathbf  C$, belongs to 
$\mathcal {LO}_{FGac}^H$ in order to prove
JEP, it is enough to observe that
$\mathbf  C \in \mathcal {LO}_{FGac}$.

To prove the last statement, given $\mathbf A$ and $\mathbf  B$,
set all the elements from $A$ to be $<$ than all the elements
from $B$ in $A \cup B$ and define the operations on $A \cup B$ 
in the unique compatible way.
\end{proof}

  Even more general conditions 
under which AP and SAP are preserved are presented
in \cite{apu}.

Recall that if $\mathcal F$
is a class of finitely generated structures
in a countable language,
a Fra\"\i ss\'e limit of  $\mathcal F$ is
a  countable  ultrahomogeneous structure  of age $\mathcal F$.
See \cite[Section 7.1]{H} for  further details.

We say that a first-order sentence $\sigma$ 
is \emph{$1$-universal} if $\sigma$ is universal
and only one variable appears in $\sigma$.
Examples of $1$-universal sentences are sentences asserting
that some unary operation is increasing, decreasing,
strictly increasing, strictly decreasing, idempotent, has
no fixed point;
that some pair of unary operations commute,  etc. 
In particular,  
there is a  $1$-universal sentence asserting that
an  order preserving unary operation 
 is  a closure operation. 

\begin{theorem} \labbel{fra}
  \begin{enumerate}   
 \item   
 Let $\mathcal K$ be 
either $\mathcal {LO}_p^H$,
 $\mathcal {LO}_r^H$  or 
  $\mathcal {PO}_F^H$
for $F$ finite, 
where $H$ is any, possibly empty,  set of conditions expressible by a $1$-universal sentence. 

If $\mathcal K^{fin}$ is the
 the class of finite members of $\mathcal K$
and $\mathcal K^{fin}$ is not empty, 
then $\mathcal K^{fin}$  
has  a Fra\"\i ss\'e limit $\mathbf M$ in $\mathcal K$.

If $\mathcal K$ is either $\mathcal {LO}_p^H$
or $\mathcal {LO}_r^H$
and 
$H$ includes the condition  $f^{m+2}(x)=f^m(x)$, for some $m$, 
 then the first-order theory $Th(\mathbf M)$ of  $\mathbf M$ is $ \omega$-categorical
and has quantifier elimination.
Moreover, $Th(\mathbf M)$ is the model-completion
of $Th(\mathcal K)$. 

\item
Let $\mathcal K$ be 
either
 $\mathcal {LO}_{p}^H$, $\mathcal {LO}_{sp}^H$,
$\mathcal {LO}_{r}^H$ or $\mathcal {LO}_{sr}^H$,  
where $H$ is any, possibly empty, set of conditions 
 expressible by a $1$-universal sentence. 

If $\mathcal K$ is nonempty and $\mathcal K^{fg}$ is 
 the class of all finitely generated members of $\mathcal K$,
then $\mathcal K^{fg}$  
has  a Fra\"\i ss\'e limit in $\mathcal K$.
 \end{enumerate}
 \end{theorem}

\begin{proof}
(1)
In  each case $\mathcal K^{fin}$
has AP and JEP, by  Lemma  \ref{inpiu}.
Obviously $\mathcal K^{fin}$ is closed under
taking  substructures,
hence the Fra\"\i ss\'e limit of $\mathcal K^{fin}$ exists
by Fra\"\i ss\'e's Theorem.  See, e.~g., \cite[Theorem 7.1.2]{H}.
The finiteness
of $F$ in   $\mathcal {PO}_F$ is necessary in order to have
only a countable number of 
nonisomorphic structures in $\mathcal K^{fin}$.

The Fra\"\i ss\'e limit belongs to $\mathcal K$ since
the limit is constructed as the union of a chain of structures in $\mathcal K^{fin}$
and $\mathcal K$ is closed under  
unions of chains.
To prove the last statement, use \cite[Theorem 7.4.1]{H},
noticing that if  $f^{m+2}(x)=f^m(x)$ holds for some $m$, then
any member of $\mathcal K$ generated by $n$
elements has cardinality $\leq (m+1)n$.    
Finally, $Th(\mathbf M)$ is model-complete
 and $\mathbf M$ is existentially closed in $\mathcal K$;
moreover, $Th(\mathcal K)$
and   $Th(\mathbf M)$ have the same universal consequences.

(2) 
is proved in  a similar way.
Just check that in each case $\mathcal K$ 
has only a countable number of nonisomorphic finitely generated members.
\end{proof}

\begin{remark} \labbel{fa}    
Fra\"\i ss\'e method does not apply
to 
the classes $\mathcal {LO}_{Fa}$
 and $\mathcal {LO}_{FGac}$, since such classes
are generally not closed under taking substructures.
The problem can be circumvented, 
since operations in $F$ and in $G$ 
are assumed to be bijective, hence we get an inessential expansion
of the language if we assume that, for every $f \in F$,
there is another operation symbol in $F$
interpreted as the inverse of $f$, and similarly
for each $g \in G$. 
Thus Theorem \ref{fra}(2)
holds for  $\mathcal {LO}^H_{Fa}$
when $F$ has two function symbols, assumed to be 
one the inverse of the other. A similar result holds for
linearly ordered set with an antiautomorphism
together with its inverse.

However, we face another problem when two or more  (anti)automorphisms
are considered, together with their inverses.
Consider $\mathbb Z \times \mathbb Z$ 
with the lexicographic order, 
let $f$ and $h$  be defined by
$f(z,w)= (z+1,w)$
and 
$h(z,w)=(z,w+n(z))$,
where $n$ is an arbitrary function
from $\mathbb Z$ to $ \{ 1, -1\} $.
If we add to the language operations
representing the inverses of    
$f$ and  $g$, then  
 $(0,0)$ generates
the whole of  $\mathbb Z \times \mathbb Z$.
Letting the function $n$ vary,
we get continuum many nonisomorphic $1$-generated structures,
hence the method in Fra\"\i ss\'e construction, as it stands, cannot be applied.

Of course, for certain sets $H$ of conditions,
it is possible that $\mathcal {LO}^H_{Fa}$
 and $\mathcal {LO}^H_{FGac}$
have only countably many models
modulo isomorphism, in which case a result analogous to 
Theorem \ref{fra}(2) holds, provided
inverses  are present in the language, as specified above.
In the case of
$\mathcal {LO}^H_{FGac}$
we also need to dispense for a constant
interpreted as the center.
\end{remark}

\begin{problem} \labbel{prob}  
Lemma \ref{inpiu}
and the proof of Theorem  \ref{fra}(1)
imply that many \emph{locally finite} theories    
of \emph{partially} ordered sets with further operations
have model completion, the simplest case
being posets with a finite number of pairwise commuting 
closure operations.
In view of the counterexamples 
in the proofs of Theorems \ref{los}(c) and \ref{rev}(b),
theories of \emph{linearly} ordered sets with  many operations generally
have not model completion.
Recall that some theory has model completion if and only if 
it has both AP and model companion.

However, it is partially an  open problem
to characterize companionable theories of linear orders
with further operations.
 \end{problem}

\begin{remark} \labbel{more}    
Theorem \ref{FJ} and Corollary \ref{oper}
can be strengthened further. We can consider many order relations
at the same time, and add conditions asserting 
that  some order is coarser than another order.
Again, conditions involving the operations can be added,
for example, conditions asserting that some operation
is increasing, or that it is idempotent.
In many cases, Theorem \ref{FJ} and Corollary \ref{oper}
apply also to binary relations which 
are not necessarily orders.
See \cite{apu} for more details. 
Moreover, Corollary \ref{oper}
holds for any number of $n$-ary operations,
with $n$ varying.   
\end{remark}

\section*{Acknowledgements} 

Work performed under the auspices of G.N.S.A.G.A. Work 
partially supported by PRIN 2012 ``Logica, Modelli e Insiemi''.
The author acknowledges the MIUR Department Project awarded to the
Department of Mathematics, University of Rome Tor Vergata, CUP
E83C18000100006.

\end{document}